\documentclass{amsproc}
\newtheorem{theorem}{Theorem}[section]
\newtheorem{lemma}[theorem]{Lemma}
\newtheorem{proposition}[theorem]{Proposition}
\newtheorem{corollary}[theorem]{Corollary}

\theoremstyle{definition}
\newtheorem{definition}[theorem]{Definition}

\newtheorem{notation}[theorem]{Notation}
\newtheorem{question}[theorem]{Question}

\theoremstyle{remark}
\newtheorem{remark}[theorem]{Remark}

\numberwithin{equation}{section}

\newcommand{\abs}[1]{\lvert#1\rvert}

\newcommand{\bd}{\noindent {\sc Proof}.\ \ }

\begin{document}

\title[The space of contact Anosov flows]
{The space of (contact) Anosov flows on  3-manifolds}

\author{Shigenori Matsumoto}

\thanks{2010 {\em Mathematics Subject Classification}. Primary 37D20.
secondary 37C40.}
\thanks{{\em Key words and phrases.} Anosov flows, contact Anosov flows,
$C^1$-open subset. }

\thanks{The author is partially supported by Grant-in-Aid for
Scientific Research (C) No.\ 20540096.}

\date{\today}

\newcommand{\AAA}{{\mathcal A}}
\newcommand{\LL}{{\mathcal L}}
\newcommand{\MCG}{{\rm MCG}}
\newcommand{\PSL}{{\rm PSL}}
\newcommand{\R}{{\mathbb R}}
\newcommand{\Z}{{\mathbb Z}}
\newcommand{\XX}{{\mathcal X}}
\newcommand{\per}{{\rm per}}
\newcommand{\N}{{\mathcal N}}

\newcommand{\PP}{{\mathbb P}}
\newcommand{\GG}{{\mathbb G}}
\newcommand{\FF}{{\mathcal F}}
\newcommand{\EE}{{\mathbb E}}
\newcommand{\BB}{{\mathbb B}}
\newcommand{\HH}{{\mathcal H}}
\newcommand{\PPP}{{\mathcal P}}
\newcommand{\UU}{{\mathcal U}}
\newcommand{\oboundary}{{\mathbb S}^1_\infty}
\newcommand{\Q}{{\mathbb Q}}

\begin{abstract}
The first half of this paper is concerned with the topology
of the space $\AAA(M)$ of (not necessarily contact)
Anosov vector fields on the unit tangent 
bundle $M$ of closed oriented hyperbolic surfaces $\Sigma$. We show that
there are countably infinite connected components of $\AAA(M)$,
each of which is not simply connected.
In the second part, we study contact Anosov flows. We show in particular
that the time changes of contact Anosov flows
form a $C^1$-open subset of the space of the
Anosov flows which leave a particular
$C^\infty$ volume form invariant, if the ambiant manifold is 
a rational homology sphere.
\end{abstract}

\maketitle

\section{Introduction}

The main purpose of this paper is to study the topology of the space 
of contact Anosov vector fields on 3-manifolds.
But going to that subject, we first consider the space $\AAA(M)$ of
(not necessarily contact) Anosov vector fields on the unit 
tangent bundle $M$ of a closed oriented hyperbolic surface $\Sigma$.

The results we obtain concerning $\AAA(M)$ are elementary and easy to show.
However the author cannot find it in the literature, which 
makes him to record these fundamental facts. Denote by $\LL(M)$
the space of nonvanishing $C^\infty$ vector fields on $M$.
There is one distiguished connected component $\LL_0(M)$ of $\LL(M)$.

\begin{theorem} \label{t1}
The space $\AAA(M)$ is contained in $\LL_0(M)$.
\end{theorem}

\begin{theorem} \label{t2}
The space $\AAA(M)$ has countably infinite connected components, each
of which is not simply connected.
\end{theorem}

After we determine the mapping class group of $M$ in Section 2, we prove these
results in Section 3.

Sections 4 and 5 are devoted to the study of contact Anosov flows.
In section 4, we determine which time change of a contact Anosov
flow is again contact Anosov. Especially we show that if
the ambiant manifold $N$ is a rational homology sphere, 
such a time change is obtained by a conjugation by an orbit 
preserving $C^\infty$ diffeomorphism.

In section 5, we study the space of contact Anosov flows.
Let $\Omega$ be a $C^\infty$ volume form on
a closed oriented manifold $N$. Denote by $\AAA_\Omega(N)$
the space of the $\Omega$-preserving Anosov vector fields.
The main result is the following.

\begin{theorem} \label{t4}
If $N$ is a rational homology sphere, the subset formed by time
changes of contact Anosov flows is $C^1$-open in $\AAA_\Omega(N)$.
\end{theorem}

In \cite{FH2}, plenty of examples of contact Anosov flows are
constructed on various manifolds including hyperbolic 3-manifolds.
Theorem \ref{t4} can also serve as producing new examples which are
$C^1$-near to classical examples.

\section{The mapping class group of $M$}

Let $\Sigma$ be a closed oriented surface of genus $\geq2$. Fix
a Riemannian metric $m_0$ of curvature $-1$. 
Let $\pi: M=T^1\Sigma\to \Sigma$ be the unit tangent bundle 
w.\ r.\ t.\ $m_0$.
The purpose of this section is to determine the mapping class group
$\MCG(M)$ of $M$, which is, by definition, the quotient of
the group of all the $C^\infty$ diffeomorphisms of $M$ by the
identity component.

Denote by $\HH$ the plane field of $M$ consisting of horizontal vectors.
The principal $S^1$ action on $M$ is denoted by $V^t$, $0\leq t\leq
2\pi$,
whose infinitesimal generator is the {\em vertical vector field} $V$.
The map $V^t$ leaves $\HH$ invarinat.
The {\em standard geodesic vector field} is a horizontal vector field
$X$
such that $\pi_*X_{(x,v)}=v$, where $v\in T^1_x(\Sigma)$ and
$x\in\Sigma$. It generates the {\em standard geodesic flow} $\{X^t\}$.

\begin{notation}
In this paper, the flow generated by a vector field $A$ is denoted by
$\{A^t\}$.
\end{notation}

The three vector fields $V$, $X$ and $Y=V^{\pi/2}_*X$ span a
Lie subalgebra,
isomorphic to the Lie algebra of $\PSL(2,\R)$.
 On the topological aspect, the following fundamental fact is a
 consequence
of the classification
of  $S^1$ bundles over surfaces. See for example \cite{O}.

\begin{proposition} \label{p1}
Any $C^\infty$ free $S^1$ action on $M$ is conjugate
to $\{V^t\}$ by a $C^\infty$ diffeomorphism isotopic to the identity.
\qed
\end{proposition}

Given $[f]\in\MCG(M)$, there is a representative $f$ which 
commutes with the $S^1$-action $\{V^t\}$. Such $f$ induces
a diffeomorphism of $\Sigma$. Thus
we get a homomorphism
$$
\pi_*: \MCG(M)\to\MCG^\diamond(\Sigma),$$
where $\MCG^\diamond(\Sigma)$ is the generalized mapping
class group of $\Sigma$ consisting of orientation preserving
or reversing classes.

Conversely, given $[g]\in \MCG^\diamond(\Sigma)$, the derivative
$dg$ yields a class $[g_*]\in\MCG(M)$, where $g_*:M\to M$ is
defined from $dg$ just by normalizing
the image vector. This yields a cross section
$$
s:\MCG^\diamond(\Sigma)\to\MCG(M).$$
Notice that $s([g])=[g_*]$ is always orientation preserving regardless
of the orientation property of $[g].$

Now let $K$ be the kernel of $\pi_*$. Any element of $K$ can be
represented
by a diffeomorphism $f$ of $M$ which preserves the fibers of the
$S^1$ action $\{V^t\}$, i.\ e.\ a diffeomorphsim
which covers the identity of $\Sigma$. Restricted to each fiber, $f$ must be
orientation preserving. For, otherwise the fixed point set of $f$
(two points set for each fiber)
would yield a multi cross section of $\pi:M\to\Sigma$,
contradicting the fact that $\pi$ is a nontrivial $S^1$ bundle.
Each class of $K$ can be 
represented by a diffeomorphism $f$ which is a rigid rotation $V^{\rho(x)}$ 
on each fiber $\pi^{-1}(x)$, where 
 $\rho:\Sigma\to S^1$ is a $C^\infty$ function.
This yields an identification
$K\cong [\Sigma,S^1]\cong H^1(\Sigma,\Z)$ and therefore we get:

\begin{proposition} \label{p2}
There is an isomorphism
\begin{equation} \label{e1}
\MCG(M)\cong H^1(\Sigma;\Z)\times\!\!\!\vert\,\, \MCG^\diamond(\Sigma).
\end{equation}
\end{proposition}

\begin{remark} \label{r1}
There is no orientation reversing homeomorphism of $M$.
\end{remark}

\section{The space $\AAA(M)$}

The vector fields $V$, $X$ and $Y$, as well as $-V$, all belong to the
same component of the space $\LL(M)$ of the nonvanishing vector fields
of $M$. Denote it by $\LL_0(M)$ and call it the {\em untwisted component}.
Notice that the components of $\LL(M)$ is in one to one correspondence
with the set $[M,S^2]$.

The differential of a diffeomorphism $f$ yields a homeomorphism
$df:\LL(M)\to \LL(M)$. 

\begin{proposition} 
For any diffeomorphism $f$ of $M$, we have $df(\LL_0)=\LL_0$.
\end{proposition}

\bd
This follows from the fact that each class of $\MCG(M)$ has a
representative which maps $V$ to a nonzero function multiple
of $V$. \qed

\smallskip
Let us denote by $\AAA(M)$ the subset of $\LL(M)$ consisiting of
Anosov vector fields.

\begin{theorem} \label{t2.1}
The space $\AAA(M)$ of the Anosov vector fields is contained in
the untwisted component $\LL_0(M)$.
\end{theorem}

\bd
In  way of showing the global structural stability theorem
for Anosov flows on the manifold $M$, E. Ghys \cite{G} proved that for any Anosov
flow
$\{A^t\}$, the weak stable foliation  can be made transverse to
the $S^1$ fibers after the conjugation by a diffeomorphism $f$. 
Each class of $\MCG(M)$ has a representative which leaves the orbit
foliation of the
$S^1$ action invariant. This implies that the conjugacy $f$ can
be chosen to be isotopic
to the identity. That is, one may assume that  the vector field $A$
which generates $\{A^t\}$ is tangent to a foliation transverse to $V$.
Then clearly the vector field $(1-s)A+sV$, $0\leq s\leq1$ is
nonvanishing, and $A$ is homotopic to $V$. \qed

\smallskip
 
Now given any negatively curved Riemannian metric
$m$ of $\Sigma$, the unit tangent bundle w.\ r.\ t.\ $m$ can be
identified with $M$ just by changing the length, and the geodesic flow
$\{A^t\}$ of
$m$ can be viewed as a flow on $M$ in the
following way. Given $p\in M$, a unit tangent vector
of $\Sigma$ w.\ r.\ t.\ $m_0$, change the length of $p$ so that the modified vector
$p'$ is a unit vector
w.\ r.\ t.\ $m$. Consider a geodesic curve $\gamma$ w.\ r.\ t.\
$m$ whose innitial velocity
vector is $p'$. Consider the vector $\gamma'(t)$ and change its
length to obtain $q\in M$. Then $A^t(p)=q$.

Let us denote by $\AAA_0(M)$ the connected component of $\AAA(M)$
which contains the standard geodesic vector field $X$.

\begin{proposition} \label{p3}
The geodesic vector field of any negatively
curved Riemannian metric on $\Sigma$ belongs to $\AAA_0(M)$.
\end{proposition}

\bd
This follows from the fact that the space of negatively curved 
Riemannian metrics is connected. \qed

\smallskip
Now for any diffeomorphism $f$ of $M$, 
we have $df(\AAA(M))=\AAA(M)$.

\begin{proposition} \label{p4}
For any element $[f]$ of $\MCG(M)$ which belongs to
 $\MCG^\diamond(\Sigma)$
in the decomposition (\ref{e1}), we have
$df(\AAA_0(M))=\AAA_0(M)$.
\end{proposition}

\bd We only need to show that for any diffeomorphism $g$ of $\Sigma$,
the induced diffeomorphism $g_*$ of $M$ carries $X$ to an element of
$\AAA_0(M)$,
i.\ e.\ $d(g_*)X\in\AAA_0(M)$. But this follows immediately from
Proposition \ref{p3}, since $d(g_*)X$ is the geodesic vector field of the Riemannian metric $(g^{-1})^*m_0$.
\qed

\smallskip

The action of $H^1(\Sigma,\Z)$ in the decomposition (\ref{e1})
on $\AAA(M)$ is quite different.
To study this we need the following lemma.

\begin{lemma} \label{l1}
Let $\{A^t\}$ be an arbitrary Anosov flow on $M$. For any
essential oriented closed curve $c$ of $\Sigma$, there
is a unique periodic orbit $\gamma$ of $\{A^t\}$ such
that $\pi(\gamma)$ is homotopic to $c$.
\end{lemma}

\bd
This is true for the standard geodesic flow $\{X^t\}$. On
the other hand any Anosov flow $\{A^t\}$ is flow 
equivalent\footnote{This means
that $h$ carries any orbit of $\{A^t\}$ onto an orbit of $\{X^t\}$
in a way to preserve the time orientation of the flows.} to $\{X^t\}$ by a homeomorphism $h$ \cite{G}.
Finally the homeomorphism $h$ can be isotoped to a diffeomorphism
$h'$
which preserve the orbit foliation of the $S^1$-action,
by an isotopy $h_t$, $0\leq t\leq 1$, where $h_0=h'$ and $h_1=h$.
Clearly the lemma holds for $\{h_0 X^t h_0^{-1}\}$.
Therefore by the continuity of the family of the topological flows, it
also holds for $\{h_1 X^t h_1^{-1}\}$. Now the latter is flow
equivalent to $\{A^t\}$,
completing the proof of the lemma. \qed

\smallskip
The next proposition shows the first half of Theorem \ref{t2}.

\begin{proposition} \label{p5}
For any nonzero element $a\in H^1(\Sigma,\Z)$, the class
$[f]$ of $\MCG(M)$ which corresponds to $a$ in (\ref{e1})
satisfies 
$df(\AAA_0(M))\cap\AAA_0(M)=\emptyset$.
\end{proposition}

\bd
We need only to show that the flow $\{f X^tf^{-1}\}$ is not
isotopic to the flow $\{X^t\}$. Choose a simple closed curve $c$
in $\Sigma$ such that $\langle a,c\rangle\neq0$.
The periodic orbit $\gamma$ in
Lemma \ref{l1} for the flow $\{X^t\}$ 
is obtained as follows. Homotope $c$
to a simple closed geodesic $l$. Then $\gamma$ is the horizontal lift
of $l$. 

Next consider the periodic orbit $\gamma'$ for the flow
$\{fX^tf^{-1}\}$. For a convenient choice of $f$ from the class, $\gamma'$ is
the image of $\gamma$ by a nontrivial Dehn twist on the torus $\pi^{-1}(l)$.
Therefore $\gamma'$ is not homotopic to $\gamma$. This shows that
the flow $\{f X^tf^{-1}\}$ is not
isotopic to the flow $\{X^t\}$. \qed

\smallskip
Let us show the last part of Theorem \ref{t2}. Let $\AAA_*$ be
an arbitrary component of $\AAA(M)$. Choose $\{A^t\}$ from
$\AAA_*$ and consider the loop $\{V^sA^tV^{-s}\}$, $0\leq s\leq 2\pi$,
in $\AAA_*$. Assume for contradiction that this loop is contractible.
Choose a periodic orbit $\gamma(t)$, $0\leq t\leq T$ of $\{A^t\}$
such that $\pi(\gamma)$ is homotopic to a simple closed
curve on $\Sigma$. Then the (possibly singular) torus
$\{V^s\gamma(t)V^{-s}\mid 0\leq t\leq T,\ 0\leq s\leq 2\pi\}$
is homotopic to an essential torus. Especially it is $\pi_1$-injective.
This contradicts that
the above loop is contractible.

\begin{remark}
We suspect that the union of $df(\AAA_0(M))$, $[f]$ from $H^1(M,\Z)$,
is the whole of $\AAA(M)$, and that $\AAA_0(M)$ is homotopy equivalent
to the circle. The analogous statement for the Anosov diffeomorphisms
on the two torus can be found in \cite{FG}. Their method is an
 application
of the thermodynamical formalism. But for flows on 3-manifolds,
it seems quite difficult to deform the Lyapunov exponent
although a potential tool is available in \cite{A}. 
\end{remark}

\section{Contact Anosov flows and their time changes}

Let $N$ be a closed oriented $C^\infty$ 3-manifold. 
\begin{definition} An Anosov flow $\{A^t\}$
(resp.\ Anosov vector field $A$) on $N$
is said to be {\em contact} if
it is the Reeb flow (resp.\ Reeb vector field) of some contact form $\tau$.
\end{definition}

If $A$ is contact Anosov, then it leaves the volume form 
$\tau\wedge d\tau$ invariant.
The $C^\infty$ plane field ${\rm Ker}(\tau)$ is invariant by $A^t$.
On the other hand the sum $E^{uu}\oplus E^{ss}$ of the strong stable
and unstable bundles is the only $A^t$-invariant
plane field transverse to $A$.  Therefore 
we have 
$${\rm Ker}(\tau)=E^{uu}\oplus E^{ss},$$
and the contact form $\tau$ is uniquely 
determined by the Anosov vector field $A$. It is known
\cite{FH1}
that if $A$ is a volume preserving Anosov flow and if $E^{uu}\oplus
E^{ss}$
is Lipschitz continuous, then 
$E^{uu}\oplus E^{ss}$ is in fact $C^\infty$, 
and the flow $A$ is contact Anosov. The contact form $\tau$
of a contact Anosov flow is shown to be tight using a result of
\cite{H}. Contact Anosov flows exhibit strong ergodic properties
\cite{L1,T1,T2}.
 The geodesic flow of a negatively curved surface is
a typical example of contact Anosov flows. In fact it
was the only known example before \cite{FH2}.

Before going to the study of time changes of contact Anosov flows,
let us recall a well known fact about the invariant volume
of an Anosov flow, which follows from the ergodicity
and the Liv\v{s}ic homological theorem
\cite{L2}.

\begin{theorem} \label{t3}
If an Anosov flow on $N$ is volume preserving, then 
the invariant volume is $C^\infty$
and unique up to a positive constant multiple.
\end{theorem}

If $A$ is an Anosov vector field and $\phi$ is a positive $C^\infty$
function, then $\phi A$ is called a {\em time change} of $A$. It is
also an Anosov vector field. If $A$ leaves the volume form $\Omega$
invariant, then $\phi A$ leaves the volume form $\phi^{-1}\Omega$
invariant.

The purpose of this section is to study what kind of time
change of a contact Anosov flow $A$ is again contact,
and the main result is Proposition \ref{p6} below. But
before going there, we need some fundamental facts.

\begin{proposition} \label{p3.1}
If $A$ is a suspension Anosov vector field, then any time change of
$A$ cannot be contact Anosov.
\end{proposition}

\bd
Since $A$ is a suspension, there is a closed 1-form $\alpha$ such
that $\alpha(A)=1$. If there is no $A$-invariant volume form,
then any time change of $A$ does not admit an invariant volume,
and it cannot be contact Anosov.
So assume $A$ admits a $C^\infty$ volume form $\Omega$.
Suppose for contradiction that $\phi A$ is contact for a contact
form $\tau$. Then by Theorem \ref{t3}, $\tau\wedge d\tau=c \phi^{-1}\Omega$
for some constant $c\neq0$. Therefore 
\begin{equation} \label{e3.1}
d\tau=\iota_{\phi A}(\tau\wedge d\tau)=\iota_{\phi A}(c \phi^{-1}\Omega)
=c\,\,\iota_A\Omega.
\end{equation}
On the other hand, since
$$
\LL_A(\alpha\wedge(\iota_A\Omega))=d\iota_A(\alpha\wedge(\iota_A\Omega))
=d\iota_A\Omega=0,$$ 
and $\alpha\wedge(\iota_A\Omega)$ is nonvanishing,
Theorem \ref{t3} shows,
\begin{equation} \label{e3.2}
\alpha\wedge(\iota_A\Omega)=c'\Omega,
\end{equation}
for some $c'\neq0$.
But $\alpha$ and $\iota_A\Omega$ are both closed, which says that
$\iota_A\Omega$ cannot be null cohomologous, contradicting (\ref{e3.1}).
\qed

\begin{proposition} \label{p3.2}
If $A$ is contact Anosov with a positive contact form, then
any time change of $A$ cannot be contact with a negative contact
form.
\end{proposition}

\bd
Assume that $A$ is contact Anosov 
with a contact form $\tau$, i.\ e.\ $\tau(A)=1$ and
$\iota_Ad\tau=0$. Also assume that there are a positive function $\phi$
and a contact form $\tau'$ such that $\tau'(\phi A)=1$,
$\iota_Ad\tau'=0$, and $\tau'\wedge d\tau'=-c\phi^{-1}\tau\wedge d\tau$
for some constant $c>0$. Then we have
$$
d\tau'=\iota_{\phi A}(\tau'\wedge d\tau')=
-c\,\,\iota_{\phi A}(\phi^{-1}\tau\wedge d\tau)
=-c\,\,\iota_A(\tau\wedge d\tau)
=-c\,\,d\tau,$$
and hence
$$\tau'=-c\tau+\omega,
$$
for some closed 1-form $\omega$.

Now for any asymptotic cycle $\Gamma$ of $A$, we have
$$
\langle\omega,\Gamma\rangle\geq\min(\tau'(A)+c\tau(A))=\min(\phi^{-1})+c>0.$$
This implies \cite{S} that $A$ has a global cross section,
contradicting Proposition \ref{p3.1}.
\qed

The following is the main result of this section.

\begin{proposition} \label{p6}
A time change $\phi A$ of a contact Anosov vector
field $A$ is again contact Anosov 
if and only if $\phi^{-1}=\omega(A)+c$ for a closed 1-form
$\omega$ and a constant $c>0$.
\end{proposition}

\bd 
Let $A$ (resp.\ $\phi A$) be a contact Anosov vector field with
the contact form $\tau$ (resp.\ $\tau'$).
Then by Proposition \ref{p3.2} we have
$\tau'=c\tau+\omega$ for some closed 1-form $\omega$ and $c>0$.
Evaluating on $\phi A$, we get
$\phi^{-1}=c+\omega(A)$.

The converse can be shown just by reversing the argument.
\qed

\smallskip
When the manifold $N$ is a rational homology sphere, the above criterion
becomes more transparant. 
Notice that there are cocompact lattices $\Gamma$ of ${\rm PSL}(2,\R)$
such that the quotient spaces $\Gamma\setminus{\rm PSL}(2,\R)$
are rational homology spheres. They all admit contact Anosov flows.
As before, let $A$ be a contact Anosov
vector field on a closed oriented 3-manifold $N$.

\begin{proposition} \label{p7}
Assume $N$ is a rational homology sphere.
A time change $B=\phi A$ is contact Anosov if and only if for some
$c>0$, the flow
$\{B^{ct}\}$ is  conjugate to $\{A^t\}$ by an
orbit preserving $C^\infty$ diffeomorphism.
\end{proposition}

\bd
For a time change $B=\phi A$ of $A$, there is a $C^\infty$ map
$a:\R\times N\to \R$ such that
\begin{equation} \label{e2}
B^{a(t,p)}(p)=A^t(p), \ \ \forall t\in\R, \ p\in N.
\end{equation}

The function $a$ is a cocycle over $\{A^t\}$, that is,
\begin{equation} \label{e3}
a(t+s,p)=a(s,A^t(p))+a(t,p).
\end{equation}
Define a function $\alpha:N\to\R$ by
$\displaystyle \alpha(p)=\frac{\partial}{\partial t}a(t,p)\vert_{t=0}$.
By (\ref{e3}), we have
\begin{equation*}
\frac{\partial}{\partial t}a(t,p)=\alpha(A^t(p)).
\end{equation*}
This implies
\begin{equation} \label{e4}
a(t,p)=\int^t_0\alpha(A^s(p))ds.
\end{equation}

Now the if part of the theorem is obvious.
So assume $B$ is also contact. Then 
by differentiating (\ref{e2}), we get $\alpha=\phi^{-1}$.
Thus Proposition \ref{p6} implies that
 $\alpha=\omega(A)+c$. Since $N$ is a 
rational homology sphere, there is a $C^\infty$ function $\psi$
such that $\omega=d\psi$, and thus
$$
\alpha=A(\psi)+c.
$$
Then (\ref{e4}) implies 
$$
a(t,p)=\psi(A^t(p))-\psi(p)+ct.
$$
Define a map $f:N\to N$ by $f(p)=B^{-\psi(p)}(p)$. Then
\begin{equation} \label{e3.3}
f(A^t(p))=B^{a(t,p)-\psi(A^t(p))}(p)=B^{ct-\psi(p)}(p)=B^{ct}(f(p)).
\end{equation}

The equation (\ref{e3.3}) implies that
the map $f$ is a $C^\infty$ diffeomorphism, showing that $\{B^{ct}\}$
is  conjugate to $\{A^t\}$ by $f$.
\qed

\section{Perturbations of a contact Anosov flow}

Let $A$ be a contact Anosov flow on a closed oriented 3-manifold
$N$, with the contact form $\tau$. Then $\Omega=\tau\wedge d\tau$
is an $A$-invariant volume form. Let us denote by $\XX_\Omega(N)$
(resp.\ $\AAA_\Omega(N)$) the space of $\Omega$-preserving vector fields 
(resp.\ $\Omega$-preserving Anosov vector fields) on $N$.
Thus $\AAA_\Omega(N)$ is a $C^1$-open subset of the linear space
$\XX_\Omega(N)$.
For any $B\in\XX_\Omega(N)$ small in the $C^1$ topology,
the flow $A+B$ is again an $\Omega$ preserving Anosov flow, i.\ e.\
$A+B\in\AAA_\Omega(N)$.
In this section we ask which $A+B$ can be a time change
of a contact Anosov flow.

Assume $\phi(A+B)$ is contact Anosov for some positive
function $\phi$, with the contact form $\tau'$.
Then by Theorem \ref{t3}, we have $\tau'\wedge d\tau'=c\phi^{-1}\Omega$
for some $c>0$.
Thus
$$
\iota_{A+B}\Omega=\iota_{\phi(A+B)}\phi^{-1}\Omega=c^{-1}
\iota_{\phi(A+B)}(\tau'\wedge d\tau')=c^{-1}d\tau',\ \mbox{ and }
$$
$$
\iota_B\Omega=\iota_{A+B}\Omega-\iota_A\Omega=c^{-1}d\tau'-d\tau,
$$
showing that $\iota_B\Omega$ is an exact 2-form.
On the other hand $B$
belongs to $\XX_\Omega(N)$ if and only if  
 $\iota_B\Omega$ is closed.
Since the correspodence $B \leftrightarrow \iota_B\Omega$ is bijective, this
show the following.

\begin{proposition} \label{p8}
The subset consisting of time changes of contact Anosov flows 
is contained  in a subspace of codimension equal to $\dim H^2(N;\R)$
in a neighbourhood of $A$ in $\AAA_\Omega(N)$.
\qed
\end{proposition}

From now on let us assume that $N$ is
a rational homology sphere and show
Theorem \ref{t4}.
Notice that 
the validity of Theorem \ref{t4} does not change if one changes $\Omega$
by a positive function multiple.
Therefore it suffices to assume that 
$A\in\AAA_\Omega(N)$ is a contact Anosov vector field 
with the contact form $\tau$ such that
$\Omega=\tau\wedge d\tau$ and to show that for any $C^1$-small
$B\in\XX_\Omega(N)$,
$A+B$ is a time change of a contact Anosov vector
field.

Now the 2-form $\iota_B\Omega$ is closed since $B$ is
$\Omega$-preserving,
and exact since $N$ is a rational homology sphere.
Choose a 1-form $\beta$ such that $d\beta=\iota_B\Omega$.
Then we have $d\tau+d\beta=\iota_{A'}\Omega$, and hence
\begin{equation} \label{e5}
\iota_{A'}(d\tau+d\beta)=0.
\end{equation}
Our goal is to show that for $C^1$-small $B$,
there is a 1-form $\tau'$ such that 
\begin{equation}\label{e6}
d\tau'=d\tau+d\beta\ \mbox{ and}\ \tau'(A')>0.
\end{equation}
 For, then
the equation (\ref{e5}), together with the fact that $\iota_{A'}\Omega$
is nonvanishing, shows that the form $\tau'$ is contact, and
the time change $\tau'(A')^{-1}A'$ 
is contact Anosov.

Since $\tau(A)=1$, we have
\begin{equation}\label{e7}
\int_\gamma\tau=\per(\gamma)\ \ \mbox{ for any periodic orbit $\gamma$
of $A$},
\end{equation}
where $\per(\gamma)$ denotes the period of $\gamma$.

Let us show that for any $\epsilon>0$, if $B$ is sufficiently $C^1$-small,
\begin{equation} \label{e8}
\int_{\gamma'}(\tau+\beta)>(1-3\epsilon)\per(\gamma')\ \ \mbox{ for any periodic orbit $\gamma'$
of $A'$}.
\end{equation}

First of all if we choose $B$ so that 
$\Vert B\Vert\,\Vert\tau\Vert<\epsilon$, then we have
$$\tau(A')=\tau(A)+\tau(B)>1-\epsilon,$$
and therefore
$$
\int_{\gamma'}\tau>(1-\epsilon)\per(\gamma')\ \ \mbox{ for any periodic orbit $\gamma'$
of $A'$}.
$$
So what we need is to show that
\begin{equation} \label{e9}
\abs{\int_{\gamma'}\beta}<2\epsilon\,\per(\gamma')\ \ \mbox{ for any periodic orbit $\gamma'$
of $A'$}.
\end{equation} 

Now the $C^1$-norm of $A'$ is bounded, regardless of the choice of $A'$
from a $C^1$-neighbourhood $\UU$ of $A$. Choose a triangulation $T$ of
$N$ by small geodesic simplices. If we choose $T$ fine enough compared with 
the above $C^1$-norm, then the orbits of $A'$ look like straight 
lines in a close-up. Thus for any periodic orbit 
$\gamma'$ of $A'$, there is a simplicial path $\gamma'_T$ and
an annulus $\mathbb A$ such that $\partial\mathbb
A=\gamma'\cup(-\gamma'_T)$
and that ${\rm Area}(\mathbb A)\leq C\,\per(\gamma')$, where
$C$ is a constant depending only on $\UU$ and $T$.

Then if $B$, and hence $d\beta=\iota_B\Omega$, is small enough, we have
\begin{equation} \label{e10}
\abs{\int_{\gamma'}\beta-\int_{\gamma'_T}\beta}
=\abs{\int_{\mathbb A}d\beta}\leq
{\rm Area}(\mathbb A)\,\Vert d\beta\Vert<\epsilon\,\per(\gamma').
\end{equation}

Denote by $\Vert\cdot\Vert_1$ the $l^1$ norm in the real coefficient
chain group
of the triangulation $T$. Then we have
\begin{equation}\label{e11}
C^{-1}\Vert \gamma'_T\Vert_1\leq\per(\gamma')\leq C\,\Vert 
\gamma'_T\Vert_1,
\end{equation}
where again $C$ depends only on $\UU$ and $T$, and is independent of 
the choice of $A'$ from $\UU$, nor of the periodic orbit $\gamma'$
of $A'$.
Now the boundary operator $\partial_2:C_2(T)\to B_1(T)$
admits a cross section $\sigma:B_1(T)\to C_2(T)$. The mapping
norm $\Vert\sigma\Vert_1$ of $\sigma$ is finite since $B_1(T)$ is finite dimensional.
Thus if $B$ is small enough, then 
\begin{equation} \label{e12}
\abs{\int_{\gamma'_T}\beta}=\abs{\int_{\sigma(\gamma'_T)}d\beta}
\leq \Vert\sigma\Vert_1\,\,\Vert \gamma'_T\Vert_1\,\,\Vert d\beta \Vert
<\epsilon \,\per(\gamma'),
\end{equation}
where the last inequality follows from (\ref{e11}).
Now (\ref{e10}) and (\ref{e12}) imply the desired inequality
(\ref{e9}). The proof of (\ref{e8}) is complete.

\smallskip
Finally let us show that (\ref{e8}) implies (\ref{e6}).
For any periodic orbit $\gamma'$ of $A'$, there is associated
an $A'$-invariant measure $\delta_{\gamma'}$ supported on $\gamma'$.
It is well known, easy to show by the specification property of Anosov
flows, that the set of measures $\delta_{\gamma'}$ is dense in
the set of the ergodic probability measures. Thus (\ref{e8})
implies that
\begin{equation} \label{e12+}
\langle\mu,(\tau+\beta)(A')\rangle\geq 1-3\epsilon
\end{equation}
for any $A'$-invariant probability measure $\mu$.

Then we have
\begin{equation} \label{e13}
t^{-1}\int^t_0(\tau+\beta)(A')\circ(A')^tdt >1-4\epsilon
\end{equation}
for any large $t$.
For, otherwise one can construct an $A'$-invariant probability
measure violating (\ref{e12+}).

If we put $$\tau'=t^{-1}\int^t_0((A')^t)^*(\tau+\beta)dt,$$
the left hand side of (\ref{e13}) coincides with $\tau'(A')$.
On the other hand, we have
$$
d\tau'=t^{-1}\int^t_0((A')^t)^*(d\tau+d\beta)dt=d\tau+d\beta,
$$
since
$$
\LL_{A'}(d\tau+d\beta)=d\,\iota_{A'}(d\tau+d\beta)=d\iota_{A'}\iota_{A'}\Omega=0.$$
This shows (\ref{e6}), as is required.
\qed

\smallskip
Theorem \ref{t4} can be generalized as follows.

\begin{corollary} \label{c1}
For an arbitrary closed 3-manifold $N$ and a $C^\infty$
volume form $\Omega$, there is a $C^1$-neighbourhood
$\UU$ of $0$ in $\XX_\Omega(N)$ such that if $A$ is a
contact Anosov vector field, and $B\in\UU$ satisfies
that $\iota_B\Omega$ is exact, then $A+B$ is a 
time change of a contact Anosov vector field.
\end{corollary}

Let $\Sigma$ be a closed oriented surface endowed with a 
Riemannian metric $m$ of varying negative curvature $K$, 
and let $\pi:M\to\Sigma$ be
the unit tangent bundle w.\ r.\ t.\ $m$. Denote the vertical vector
field by $V$, the geodesic vector field by $X$, and $Y=V^{\pi/2}_*X$.
They satify:
\begin{equation} \label{e5.1}
[V,X]=Y,\ \ [V,Y]=-X,\ \ [X,Y]=K\circ\pi\, V.
\end{equation}
The 1-forms $\xi$, $\eta$ and $\theta$ dual to $X$, $Y$ and $V$ satisfy:
\begin{equation} \label{e5.2}
d\xi=\theta\wedge\eta,\ \ d\eta=-\theta\wedge\xi,\ \ d\theta=-K\circ\pi\,
\xi\wedge\eta.
\end{equation}
The volume form $\Omega=\xi\wedge\eta\wedge\theta$ is left invariant
by the three vector fields $V$, $X$ and $Y$.

G. P. Paternain \cite{P} considers what is called the {\em 
magnetic vector field}
 $$A_\lambda=X+\lambda V$$ for a costant $\lambda$, 
and shows the following.

\begin{theorem} \label{t6} {\bf (G. P. Paternain)}
\ For $\abs{\lambda}$ small, the vector field
$A_\lambda$ is not contact, unless $K$ is constant.
\end{theorem}

Let us consider 
more generally the vector field 
$$A_\phi=X+\phi\circ\pi\, V$$ for
a $C^\infty$ function $\phi:\Sigma\to\R$. 

Now we have
\begin{equation}\label{e5.3}
\LL_{A_\phi}\Omega=d\,\iota_{A_\phi}\Omega=
d(\phi\circ\pi\,\xi\wedge\eta)=0,
\end{equation}
where the last equality follows from
$
V(\phi\circ\pi)=0$. Thus the vector field $A_\phi$ leaves
the volume form $\Omega$ invariant, and it is Anosov for
$C^1$-small $\phi$. Applying Corollary \ref{c1}, we get:

\begin{proposition} \label{last}
For any negatively curved metric $m$
and a $C^1$-small function $\phi:\Sigma\to\R$, the vector field $A_\phi$ is a time change
of a contact Anosov vector field.
\end{proposition}
\bd What we need to show is that the closed 1-form 
$\phi\circ\pi\,\xi\wedge\eta$ is exact, which is
an easy consequence of the fact that $H_2(M,\Z)$ is
generated by vertical tori and that
$$\iota_V(\phi\circ\pi\,\xi\wedge\eta)=0.$$ 
\qed

\smallskip
The contact forms which appear in Proposition \ref{last} are
$C^1$-perturbations of the contact form $\xi$ and is positive.
On the other hand, the connection form $\theta$ is negative and
tight by a result of \cite{H}. Compare Remark \ref{r1}.

\begin{question}
Is there a contact Anosov flow on $M$ whose contact form is $\theta$?
\end{question}


\begin{thebibliography}{99}
\bibitem[A]{A} M. Asaoka, \it On invariant volumes of codimension-one
	   Anosov flows and the Verjovsky conjecture, \rm Invent.\
	   Math.\ {\bf 174}(2008), 435-462.

\bibitem[FG]{FG} F. T. Farrell and A. Gogolev, \it The space of Anosov
diffeomorphisms, \rm ArXiv:1201.3595.

\bibitem[FH1]{FH1} P. Foulon and B. Hasselblatt, \it Zygmund strong
foliations, \rm Israel J. Math,\ {\bf 138}(2003), 157-169.

\bibitem[FH2]{FH2} P. Foulon and B. Hasselblatt, \it Contact Anosov
flows on hyperbolic 3-manifolds, \rm Preprints, Tufts University.

\bibitem[G]{G} E. Ghys, \it Flots d'Anosov sur les 3-vari\'et\'es
	   fibr\'ees en cercles, \rm Ergo.\ Th.\ Dyn.\ Sys.\ {\bf
 4}(1984), 67-80.

\bibitem[H]{H} H. Hofer, \it Pseudoholomorphic curves in
	   symplectizations with applications to the Weinstein conjecture
	   in dimension three, \rm Invent.\ Math.\ {\bf 114}(1993), 515-563.

\bibitem[L1]{L1} C. Liverani, \it On contact Anosov flows, \rm Ann.\
	   Math.\ {\bf 159}(2004), 1275-1312.

\bibitem[L2]{L2} A. Liv\v{s}ic, \it Certain properties of the homology of
	   $Y$-systems, \rm Akademiya Nauk SSSR, Matematicheskie Zametki
	   {\bf 10}(1971), 555-564; \it Cohomology of dynamical systems,
	   \rm Izvestiya Akademii Nauk SSSR, Seriya Matematicheskaya
	   {\bf 36}(1972), 1296-1320.


\bibitem[O]{O} P. Orlik, ``Seifert manifolds,'' Lecture Notes in
	Mathematics, 291, Springer Verlag 1972.

\bibitem[P]{P} G.  P. Paternain, \it The longitudinal KAM-cocycle
of a magnetic flows, \rm Math.\ Proc.\ Cambridge Philos.\ Soc.\ {\bf
	   139}(2005), 307-316.

\bibitem[S]{S} S. Schwartzman, \it Asymptotic cycles, \rm Ann.\ Math.\
{\bf 66}(1957), 270-284.

\bibitem[T1]{T1} M. Tsujii, \it Quasi-compactness of transfer operators
	   for contact Anosov flow, \rm Nonlinearity {\bf 23}(2010),
	   1495-1545.
\bibitem[T2]{T2} M. Tsujii, \it Contact Anosov flows and the FBI
	   transform, \rm Ergod.\ Th.\ Dyn.\ Sys.\ {\bf 32}(2012), no.\ 6.

\end{thebibliography}
\end{document}